\documentclass{article}
\usepackage{amssymb, verbatim}
\usepackage{hyperref}
\usepackage{amsmath}
\usepackage{tikz-cd}

\def\Ric{\mathop{\rm Ric}}
\def\cRic{\mathop{\rm R\i\makebox[0pt]{\raisebox{5pt}{\tiny$\circ$\;\,}}c}}

\def\dist{\mathop{\rm dist}}

\def\Riem{\mathop{\rm Rm}}

\def\Vol{\mathop{\rm Vol}}

\def\RRR{\mathop{\mathcal{R}}}

\def\CP{\mathop{\mathbb{CP}}}

\def\be{\begin{eqnarray}}
\def\ee{\end{eqnarray}}
\def\beg{\begin{eqnarray*}}
\def\ees{\end{eqnarray*}}

\def\Xint#1{\mathchoice
{\XXint\displaystyle\textstyle{#1}}%
{\XXint\textstyle\scriptstyle{#1}}%
{\XXint\scriptstyle\scriptscriptstyle{#1}}%
{\XXint\scriptscriptstyle\scriptscriptstyle{#1}}%
\!\int}
\def\XXint#1#2#3{{\setbox0=\hbox{$#1{#2#3}{\int}$ }
\vcenter{\hbox{$#2#3$ }}\kern-.5\wd0}}

\newcommand{\qed}{\hfill$\Box$}
\newtheorem{theorem}{Theorem}[section]
\newtheorem{proposition}[theorem]{Proposition}
\newtheorem{lemma}[theorem]{Lemma}

\newtheorem{corollary}[theorem]{Corollary}

\date{\small\it September 12, 2013}
\author{Brian Weber}
\title{Curvature Estimates for Critical 4-Manifolds with a Lower Ricci Curvature Bound}

\begin{document}

\maketitle

\section{Introduction}

This paper follows \cite{Web3}, where certain Harnack-style inequalities and some consequences were proved on low energy regions of manifolds with an elliptic system.
Those results were of collapsing type, meaning independent of volume ratios, isoperimetric inequalities, or Sobolev constants, but were restricted to the local scale, or definite multiples thereof.
In this paper, we use those results and some additional techniques to prove elliptic-type estimates on fixed scales.

The elliptic systems we consider are the same as those in \cite{Web3}.
Our Riemannian manifolds $(N^4,g)$ are all 4-dimensional, and are either Einstein, CSC (constant scalar curvature) half-conformally flat, CSC Bach-flat, extremal K\"ahler, or have some other elliptic system in the curvature.
In all cases we assume a uniform lower bound on Ricci curvature: $\Ric\ge-\Lambda^2$.
In the extremal K\"ahler case, we also require a bound on the gradient of scalar curvature: $|\nabla{R}|\le\Lambda^3$.
In the case of any other type of metric with an elliptic system, we actually require two-sided bounds on Ricci curvature: $|\Ric|\le\Lambda^2$, or the results of \cite{Web3} will be unavailable.
These assumptions apply throughout the paper, except in Section \ref{SubSecGromovStylePacking} where a Gromov-style covering lemma is proved for Riemannian metrics, without apriori assumptions.

Before stating our results, we require a few definitions.
The first is the familiar curvature radius with a cutoff:
\be
r^s_{\mathcal{R}}(p)\;=\;\sup\,\left\{\;r\,\in\,(0,\,s)\;\big|\;|\Riem|<r^{-2}\;{\rm on}\;B(p,\,r)\;\right\}.
\ee
This is also termed the $s$-local curvature radius; compare (1.7) and (2.2) of \cite{CT}.
We use $r_{\mathcal{R}}$ for $r_{\mathcal{R}}^\infty$.
If $\Omega\subset{N}$ is a subset of a Riemannian manifold, define its $s$-thickening to be
\be
\Omega^{(s)}\;=\;\left\{\,x\in{N}\;\big|\;\dist(x,\,\Omega)\,<\,s\;\right\}.
\ee
To improve the sharpness of our first proposition, we shall also define
\be
\Omega^{\mathcal{R},s}\;=\;\left\{\,x\in{B}\left(p,\,r^s_{\mathcal{R}}(p)\right)\;\text{for some}\;p\in\Omega\;\right\}.
\ee
Of course $\Omega^{\mathcal{R},s}\subset\Omega^{(s)}$.
We use $\Omega^{\mathcal{R}}$ for $\Omega^{\mathcal{R},\infty}$.
\begin{proposition}[Integration Lemma] \label{PropIntegrationLemma}
Given $K<\infty$ and $k>0$, there exist numbers $\epsilon_0=\epsilon(K)>0$, $C=C(K,k)<\infty$, and a universal constant $m<\infty$, so that $s<{K}\Lambda^{-1}$ and $\int_{\Omega^{(2s)}}|\Riem|^2\;\le\;\epsilon_0$ together imply
\be
\int_{\Omega}\left(r^s_{\mathcal{R}}\right)^{-k}
\;\le\;m\,{s}^{-k}\,\Vol\,\Omega^{(s)}
\;+\;{C}\;\int_{\Omega^{\mathcal{R},s}}|\Riem|^{\frac{k}{2}}.
\ee
Further, if $\int_{\Omega^{2\mathcal{R}}}|\Riem|^2<\epsilon_0$ and either $\Ric\ge0$ or else $r_{\mathcal{R}}(p)<K\Lambda^{-1}$ for all $p\in\Omega$, we have
\be
\int_{\Omega}\left(r_{\mathcal{R}}\right)^{-k}
\;\le\;{C}\;\int_{\Omega^{\mathcal{R}}}|\Riem|^{\frac{k}{2}}
\ee
where $C$ does not depend on $K$ in the case $\Ric\ge0$.
\end{proposition}

\begin{theorem}[Epsilon Regularity] \label{ThmEpsReg}
Pick $K<\infty$.
There exists numbers $\epsilon_0>0$, $C<\infty$ so that
\be
\int_{B(p,r)}|\Riem|^2\;\le\;\epsilon_0
\ee
implies
\begin{itemize}
\item[{\it{i}})] If $r\ge{K}\Lambda^{-1}$, either $\sup_{B(p,\frac12r)}|\Riem|<C\Lambda^2$, or $\Xint{\;-}_{B(p,\Lambda^{-1})}\left(\frac13{R}^2+4\left|W^+\right|^2\right)>\Lambda^4$ and 
$$|\Riem(q)|^2<C\Xint{\;-}_{B(q,\Lambda^{-1})}\left(\frac13{R}^2+4\left|W^+\right|^2\right)$$
for all $q\in{B}(p,\frac12r)$.
Here $C$ and $\epsilon_0$ are universal.
\item[{\it{ii}})] If $r<K\Lambda^{-1}$, either $\sup_{B(p,\frac12r)}|\Riem|<Cr^{-2}$, or $\Xint{\;-}_{B(r)}\left(\frac13{R}^2+4\left|W^+\right|^2\right)>Cr^{-4}$ and 
$$\sup_{B(p,\frac12r)}|\Riem|^2<C\Xint{\;-}_{B(p,r)}\left(\frac13{R}^2+4\left|W^+\right|^2\right)$$
for all $q\in{B}(p,\frac12r)$, where $C=C(K)$, $\epsilon_0=\epsilon_0(K)$.
\item[{\it{iii}})] Again if $r<K\Lambda^{-1}$ then $C=C(K)$, $\epsilon_0=\epsilon_0(K)$ and
$$
\sup_{B(p,\frac12r)}|\Riem|^2\;<\;C\,\Xint{\;-}_{B(p,r)}|\Riem|^2.
$$
\end{itemize}
\end{theorem}

{\bf Remark}.
This should be compared to the standard $\epsilon$-regularity result, discussed in Lemma \ref{LemmaStdEpsReg}, which relies on estimating Sobolev constants, and to Theorem 0.8 of \cite{CT}.

{\bf Remark}.
When the metric is half-conformally flat or extremal K\"ahler, the expression $\frac13{R}^2+4\left|W^+\right|^2$ reduces to $\frac13R^2$ or $R^2$.
In the extremal K\"ahler case, one naturally has a priori controls on $\int{R}^2$.

\begin{corollary}[Epsilon-regularity for the curvature radius] \label{CorHarnack}
Given $K<\infty$ and $k>0$, there exist numbers $\epsilon_0=\epsilon_0(K)>0$ and $C=C(K)<\infty$ so that $r<K\Lambda^{-1}$ and $\int_{B(q,\,r)}|\Riem|^2\;\le\;\epsilon_0$ implies
\be
\sup_{B(p,\frac12r)}r_{\mathcal{R}}(q)^{-4}\;\le\;C\Xint{\;-}_{B(q,\,r)}\left(r_{\mathcal{R}}\right)^{-4}. \label{IneqHarnIIa}
\ee
\end{corollary}
Of course Corollary \ref{CorHarnack} can be used with the integration lemma or with Corollary \ref{CorEllipticEsts} to obtain additional inequalities.

{\bf Remark}.
Constant scalar curvature is not required in the half-conformally flat or Bach flat cases, just a bound on its gradient: $|\nabla{R}|<\Lambda^3$.

{\bf Remark}. 
In the extremal K\"ahler and CSC half-conformally flat cases, we conjecture that Corollary \ref{CorHarnack} and ({\it{iii}}) of Theorem \ref{ThmEpsReg} hold without the assumption that $\Ric\ge-\Lambda^2$.

\section{Prior Results}

For convenience, we briefly present the prior results necessary to our work.
We will, however, present only a few definitions surrounding N-structures, despite their crucial importance; we refer the reader to Section 2.2 of \cite{Web3} for a full description of the conventions used in both that paper and this.

\subsection{Harnack-style results on the curvature radius}

\begin{theorem}[Harnack inequality I] \label{ThmHarnackI}
Given $K<\infty$, there exists $\epsilon_0=\epsilon_0(\Lambda,K)>0$ and $\delta_0=\delta_0(\Lambda,K)>0$, so that $r_{\mathcal{R}}(p)<K\Lambda^{-1}$ and $\int_{B(p,\,2Kr_{\mathcal{R}}(p))}|\Riem|^2\;\le\;\epsilon_0$ together imply $r_{\mathcal{R}}\;\ge\;\delta_0\,r_{\mathcal{R}}(p)$
on $B(q,Kr_{\mathcal{R}}(p))$.
\end{theorem}

\begin{theorem}[Harnack inequality II] \label{ThmHarnackII}
Given $K<\infty$, $k>0$, and $\mu\in(0,1]$, there exist numbers $\epsilon_0=\epsilon_0(K)>0$ and $C=C(\mu,k,K)<\infty$ so that $r_{\mathcal{R}}(p)<K\Lambda^{-1}$ and $\int_{B(q,\,2{r}_{\mathcal{R}}(p))}|\Riem|^2\;\le\;\epsilon_0$ together imply
\be
r_{\mathcal{R}}(p)^{-2k}\;\le\;\frac{C}{\Vol\,B(p,\,\mu{r}_{\mathcal{R}}(p))}\int_{B(p,\,\mu{r}_{\mathcal{R}}(p))}|\Riem|^{k}. \label{IneqHarnIIa}
\ee
\end{theorem}
Simply using the estimate first on $\left(r_{\mathcal{R}}\right)^{-2}$, we also obtain
\be
r_{\mathcal{R}}(p)^{-2k}\;\le\;C\left(\frac{1}{\Vol\,B(p,\,\mu{r}_{\mathcal{R}}(p))}\int_{B(p,\,\mu{r}_{\mathcal{R}}(p))}|\Riem|^2\right)^{\frac{k}{2}}. \label{IneqHarnIIb}
\ee
\begin{corollary}[Elliptic estimates for the curvature radius] \label{CorEllipticEsts}
Given $K<\infty$ and $l\in\mathbb{N}$, there is an $\epsilon_0=\epsilon_0(\Lambda,K)>0$ and $C=C(\Lambda,l,K)$ so that if $r_{\mathcal{R}}(p)<K\Lambda^{-1}$ and
\be
\int_{B(p,\,2Kr_{\mathcal{R}}(p))}|\Riem|^2\;\le\;\epsilon_0
\ee
then
\be
\sup_{B(p,Kr^s_{\mathcal{R}}(p))}\left|\nabla^lr_{\mathcal{R}}\right|\;\le\;C\,\left(r^s_{\mathcal{R}}(p)\right)^{1-l}.
\ee
\end{corollary}
When $l\ge1$, $\mu\in(0,1]$ and $r_{\mathcal{R}}(p)<K\Lambda^{-1}$ and $\int_{B(p,2r_{\mathcal{R}}(p))}|\Riem|^2<\epsilon_0$, then using (\ref{IneqHarnIIa}) and (\ref{IneqHarnIIb}) we have
\be
\begin{aligned}
\sup_{B(p,r_{\mathcal{R}}(p))}\left|\nabla^lr_{\mathcal{R}}\right|
&\;\le\;\frac{C}{\Vol\,B(p,\,\mu{r}_{\mathcal{R}}(p))}\int_{B(p,\mu{r}_{\mathcal{R}}(p))}|\Riem|^{\frac{l-1}{2}} \\
\sup_{B(p,r_{\mathcal{R}}(p))}\left|\nabla^lr_{\mathcal{R}}\right|
&\;\le\;C\left(\frac{1}{\Vol\,B(p,\,\mu{r}_{\mathcal{R}}(p))}\int_{B(p,\mu{r}_{\mathcal{R}}(p))}|\Riem|^2\right)^{\frac{l-1}{4}}.
\end{aligned}
\ee

\subsection{Epsilon-regularity, collapsing, and lower Ricci bounds}

\begin{lemma}[Standard epsilon regularity] \label{LemmaStdEpsReg}
Given $K<\infty$ there exists numbers $\epsilon_0=\epsilon_0(K)>0$, $C=C(K)<\infty$ so that $r\le{K}\Lambda^{-1}$ and
\be
\Xint{\;-}_{B(p,r)}|\Riem|^2\;\le\;\epsilon_0\,r^{-4} \label{IneqEpsRegAssumption}
\ee
imply
\be
\sup_{B(p,r/2)}|\Riem|^2\;\le\;C\Xint{\;-}_{B(p,r)}|\Riem|^2. \label{IneqEpsRegConclusion}
\ee
\end{lemma}

{\bf Remark}.
We use Lemma \ref{LemmaStdEpsReg} iin proving the more general epsilon-regularity of Theorem \ref{ThmEpsReg}.

{\bf Remark}.
Lemma \ref{LemmaStdEpsReg} is normally proved by estimating the Sobolev constant, which is possible by using the almost-monotonicity of volume ratios when $r<K\Lambda^{-1}$, and then running through the Moser iteration process.
But it is can be proven this Lemma using only volume comparison, the Harnack inequality (Theorem \ref{ThmHarnackI}), and the integration lemma that we prove below, which is itself a consequence of Theorem \ref{ThmHarnackI}; this is done in Section \ref{SubSectionPfiii} below.
Now the proof of Theorem \ref{ThmHarnackI} uses regularity theory in an essential way, in that $C^2$ convergence of the metric is required when forming blow-up models.
But this does not involve any estimation of Sobolev constants.
Additionally, the proof in \cite{Web3} uses Lemma \ref{LemmaStdEpsReg} explicitly in one stage, but non-crucially.

An argument from \cite{CT} shows the use of this theorem in collapsing theory, as a way of enforcing collapse on curvature scales.
The significance of this is explained in Section \ref{SubSecCollapsingStructs}.
\begin{lemma}[Low energy collapse] \label{LemmaChTiForcedCollpasing}
Given $K<\infty$, $\tau>0$, there is an $\epsilon=\epsilon(\tau,K)>0$ so that $r_{\mathcal{R}}(p)\le{K}\Lambda^{-1}$ and $\int_{B(p,2r_{\mathcal{R}}(p))}|\Riem|^2\;\le\;\epsilon$ imply $\Vol\,B(p,\,r_{\mathcal{R}}(p))\;\le\;\tau\cdot{r}_{\mathcal{R}}(p)^{4}$.
\end{lemma}
\underline{\sl Pf}.
Assuming that $r_{\mathcal{R}}(p)^{-4}\,Vol\,B(p,r_{\mathcal{R}}(p))>\tau$, then choosing $r=2r_{\mathcal{R}}(p)$ and using volume comparison, we have $Vol\,B(p,r)>C\tau{r}^{-4}$ for some $C=C(K)$.
Then choosing $\epsilon$ small enough we obtain (\ref{IneqEpsRegAssumption}), so the conclusion (\ref{IneqEpsRegConclusion}) holds.
Then
\be
\begin{aligned}
r_{\mathcal{R}}(p)^{-2}
&\;=\;\sup_{B(p,r/2)}|\Riem|
\;\le\;C\left(\frac{1}{Vol\,B(p,r)}\int_{B(p,r)}|\Riem|^2\right)^{\frac12} \\
&\;\le\;C\left(\frac{\epsilon}{Vol\,B(p,2r_{\mathcal{R}}(p))}\right)^\frac12
\end{aligned}
\ee
so $r_{\mathcal{R}}(p)^{-4}Vol\,B(p,r_{\mathcal{R}}(p))<C^2\epsilon$.
Possibly choosing $\epsilon$ smaller, we again have 
$$
r_{\mathcal{R}}(p)^{-4}Vol\,B(p,r_{\mathcal{R}}(p))\le\tau.
$$
\qed

\subsection{Collapsing structures} \label{SubSecCollapsingStructs}

The F- and N-structures of Cheeger-Gromov \cite{CG1} \cite{CG2} and Cheeger-Gromov-Fukaya \cite{CFG} will be important in what follows.
For brevity, we will not define all useful aspects of these structures, but only describe the aspects that will be of use to us in Section \ref{SubSecTransgression}.
Our definition is similar to the one used in \cite{CR}; one can find a full description of our conventions in Section 2 of \cite{Web3}.

An N-structure $\mathfrak{N}$ is a triple $(\Omega,\,\mathcal{N},\,\iota)$ where $\Omega$ is a domain in a differentiable manifold, $\mathcal{N}$ is a sheaf of nilpotent Lie algebras on $\Omega$, and $\iota:\mathcal{N}\rightarrow\mathcal{X}(\Omega)$ (called the action) is a sheaf monomorphism from $\mathcal{N}$ into the Lie algebra sheaf $\mathcal{X}(\Omega)$ of differentiable vector fields on $\Omega$, so that a collection of sub-structures $\mathcal{A}=\{(\mathcal{N}_i,\Omega_i,\iota_i)\}_i$ exists that satisfies the three conditions below.
In what follows, if $p\in\Omega_i$ its $\mathcal{N}_i$-stalk will be denoted $\mathcal{N}_{i,p}$ and its $\mathcal{N}$-stalk will be denoted $\mathcal{N}_p$.
\begin{itemize}
\item[{\it{i}})] (Completeness of the cover) The collection $\{\Omega_i\}$ is a locally finite cover of $\Omega$, and given $p\in\Omega$ there is at least one $\Omega_i$ so that $\mathcal{N}_{i,p}=\mathcal{N}_p$.
\item[{\it{ii}})] (Uniformity of the action) The lifted sheaf $\widetilde{\mathcal{N}_i}$ over the universal cover $\widetilde{\Omega_i}\rightarrow{\Omega_i}$ is a constant sheaf (each stalk is canonically isomorphic to the Lie algebra of global sections $\widetilde{\mathcal{N}_i}(\widetilde{\Omega_i})$).
\item[{\it{iii}})] (Integrability of the action) Given $\Omega_i$, there is a connected, simply-connected nilpotent Lie group $G_i$ so that for any choice of $\widetilde{\iota_i}$ there is an action of $G_i$ on $\widetilde{\Omega_i}$ whose derived action is equal to the image of the Lie algebra of sections $\widetilde{\mathcal{N}_i}(\widetilde{\Omega_i})$ under $\widetilde{\iota_i}$.
\end{itemize}
An N-structure is an called an F-structure if the associated sheaf $\mathcal{N}$ is abelian.
There is a difference between our definition of F-structures and the common definition: we do not require that a {\it torus} acts on a {\it finite normal} cover of $\Omega_i$, but that $\mathbb{R}^n$ acts on its {\it universal} cover.
This is a slight convenience in what follows, as we will frequently pass to universal covers and want to refer to the lifted structures as F-structures.
The main disadvantage of using this definition is that a manifold with an F-structure need not have vanishing Euler number; for instant $\mathbb{R}^n$ with a translational vector field provides an F-structure.

\begin{lemma}[Global integrability for N-structures \cite{Web3}] \label{LemmaUniqueGlobalExtensions}
If $\mathfrak{N}$ is any N-structure on a domain $\Omega$ with an invariant metric, and if $\Omega$ is simply connected, then any element $\mathfrak{b}\in\mathcal{N}_p$ of the stalk at any point $p\in\Omega$ extends uniquely to a Killing field $V$ on $\Omega$.
\end{lemma}

\begin{theorem} \label{ThmCollapseBigThm}
There exists a $\tau>0$ so that if for all $p\in\Omega^{(s)}$ we have 
$$
\left(r_{\mathcal{R}}^s(p)\right)^{-4}Vol\,B(p,r_{\mathcal{R}}^s(p))<\tau,
$$
then there exists an N-structure $\mathfrak{N}$ on a neighborhood $\widetilde\Omega$ of $\Omega$ with the following properties:
\begin{itemize}
\item[{\it{a}})] {\rm (Saturation \cite{CFG})} $\widetilde\Omega$ is saturated with respect to the orbits of $\mathfrak{N}$, and $\widetilde\Omega\subset\Omega^{(s)}$.
\item[{\it{b}})] {\rm (Polarizability \cite{Rong})} $\widetilde\Omega$ is polarizable, with polarized atlas $(\Omega_i,\mathfrak{N}_i,\iota_i)$.
\item[{\it{c}})] {\rm (Regularity of the Atlas \cite{CG1})} On overlaps the stalks of the subsheaves $\mathcal{N}_i$ can be ordered by strict inclusion: if $p\in\Omega$ and $\{\Omega_{i_j}\}_{j=1}^k$ are the $\Omega_i$ that contain $p$, then after re-arrangement the stalks satisfy $\mathcal{N}_p=\mathcal{N}_{i_1,p}\subset\dots\subset\mathcal{N}_{i_k,p}$
\item[{\it{d}})] {\rm ($C$-regularity on the curvature scale \cite{CFG} \cite{CT})} The multiplicity of the cover $\{\Omega_i\}$ is bounded by $C(n)$, and the orbits of each substructure $(\Omega_i,\mathfrak{N}_i,\iota_i)$ have second fundamental forms at $p$ bounded by $C(n)r_{\mathcal{R}}^s(p)^{-1}$.
\item[{\it{e}})] {\rm (Almost-invariance \cite{CFG})} If $p\in\widetilde\Omega$ and the metric scaled so $r_{\mathcal{R}}(p)=1$, then the metric is close, depending on $\tau$, in the $C^{k,\alpha}$-sense to an invariant metric.
\item[{\it{f}})] {\rm (Completeness \cite{CFG} \cite{NT1})} If any of the $\widetilde\Omega_i$ is given its invariant metric, then passing to the orbit space is a projection onto an orbifold whose injectivity radius at any point is bounded from below in terms of the point's curvature radius.
\end{itemize}
Further, the center of an N-structure of positive rank is an F-structure of positive rank.
\end{theorem}


\section{Proof of the Integration Lemma} \label{SectionPfIntLemma}

\subsection{A Gromov-Style Covering Lemma} \label{SubSecGromovStylePacking}

For the following two lemmata we forget the lower Ricci bound and consider $n$-dimensional Riemannian manifolds $M^n$ with no additional conditions.
Let $\Omega\subseteq{M}^n$ be any domain.
We say a collection of points $\{p_i\}_{i=1}^N\subset\Omega$ (where possibly $N=\infty$) is a maximally $\frac1kr_{\RRR}$-separated subset if 
\be
\dist(p_i,p_j) \;\ge\; \frac1k\max\{{r}_{\RRR}(p_i),\,r_{\RRR}(p_j)\},
\ee
and if $p\in{K}$ then some $j\in\{1,\dots,N\}$ exists so that
\be
\dist(p,\,p_j) \;\le\; \frac1k\max\{{r}_{\RRR}(p),\,r_{\RRR}(p_j)\}.
\ee
In other words, for any $i$ the ball $B(p_i,\,\frac1kr_{\mathcal{R}}(p_i))$ contains no $p_j$ except $p_i$, and if $p\in\Omega$, then either $B(p,r_{\mathcal{R}}(p))$ contains one of the $p_i$ or there is a $p_i$ so that $B(p_i,r_{\mathcal{R}}(p_i))$ contains $p$.
Note that the $B(p_i,\frac{1}{2k}r_{\mathcal{R}}(p_i))$ are disjoint.

\begin{lemma}[Gromov-Style Covering Lemma] \label{LemmaGromovCoverI}
Let $\Omega\subseteq{N}^n$ be a domain in a Riemannian $n$-manifold, and let $\{p_i\}$ be a maximally $\frac1k{r}_{\RRR}$-separated subset of $\Omega$ where $k>1$.
Then if $l>\frac{k}{k-1}$
\be
\Omega \;\subseteq\; \bigcup_{i=1}^N\;B\left(p_i,\,\frac{l}{k}r_{\RRR}(p_i)\right).
\ee
If in addition $l\le\frac17k$ (so also $k\ge8$), then the multiplicity of this cover is uniformly bounded by $C(n)\cdot{k}^n$.
\end{lemma}
{\it Proof}. \; Assume $p\in\Omega\setminus\bigcup_iB_{p_i}(\frac{l}{k}r_{\RRR}(p_i))$.
There is a $j$ so that
\be
&&\dist(p,\,p_j) \;\le\; \frac1k\max\{{r}_{\RRR}(p),\,r_{\RRR}(p_j)\}
\ee
but so $\frac{l}{k}r_{\RRR}(p_j)\;\le\;\dist(p,\,p_j)$.
Then
\be
&&\frac{l}{k}r_{\RRR}(p_j)\;\le\;\dist(p,\,p_j) \;\le\; \frac1k\max\{{r}_{\RRR}(p),\,r_{\RRR}(p_j)\} \label{IneqCoverBasic}
\ee
and with $l>1$ this implies $\max\{r_{\RRR}(p),r_{\RRR}(p_j)\}=r_{\RRR}(p)$.
Because $r_{\mathcal{R}}$ is Lipschitz with Lipschitz constant $1$ or less, 
\be
\begin{aligned}
r_{\RRR}(p_j)
&\;\ge\;r_{\RRR}(p)\,-\,\dist(p,\,p_j) \\
&\;\ge\; r_{\RRR}(p)\,-\,\frac1kr_{\RRR}(p) \;=\;\frac{k-1}{k}\,r_{\RRR}(p).
\end{aligned}
\ee
However, from (\ref{IneqCoverBasic}) we have $\frac{l}{k}r_{\RRR}(p_j)\le\frac1kr_{\RRR}(p)$ so that
\be
\frac{l}{k}r_{\RRR}(p_j) &\le& \frac{1}{k-1}r_{\RRR}(p_j)
\ee
meaning $l\le\frac{k}{k-1},$ which we had assumed to be false.
This proves the first statement.

For the second assertion, consider a point $p'\in\Omega$, and let $\{p'_j\}_{j=1}^m$ be the set of points so that $p'\in{B}(p'_k,\,\frac{l}{k}r_{\RRR}(p'_k))$.
Let $p''\in\{p'_j\}_j$ be a point with $r_{\RRR}(p'')=\max_j\{r_{\RRR}(p'_j)\}$.
Then
\be
\bigcup_j{B}\left(p'_j,\,\frac{l}{k}r_{\RRR}(p'_j)\right) \;\;\subset\;\; B\left(p'',\,3\frac{l}{k}r_{\RRR}(p'')\right)
\ee
Therefore, because the $B(p'_j,\,\frac{1}{2k}r_{\RRR}(p'_j))$ are disjoint, we have
\be
\begin{aligned}
&\sum_i\Vol{B}\left(p'_j,\,\frac{1}{2k}r_{\RRR}(p'_j)\right) \;\le\; \Vol B\left(p'',\,3\frac{l}{k}r_{\RRR}(p'')\right) \\
&m\cdot\min_j\Vol{B}\left(p'_j,\,\frac{1}{2k}{r}_{\RRR}(p'_j)\right) \;\le\; \Vol B\left(p',\,3\frac{l}{k}r_{\RRR}(p'')\right).
\end{aligned} \label{IneqMultBound}
\ee
Again using $Lip(r_{\mathcal{R}})\le1$, and since $\dist(p'',p'_j)\le\frac{l}{k}(r_{\RRR}(p'')+r_{\RRR}(p'_j))$, we have 
\be
r_{\RRR}(p'_j)>r_{\RRR}(p'')-\dist(p'',p'_j) 
\;\ge\; r_{\RRR}(p'')-\frac{l}{k}(r_{\RRR}(p'')+r_{\RRR}(p'_j)) 
\ee
which implies a sandwich: for any of the $j\in\{1,\dots,m\}$ we have
\be
\frac{k-l}{k+l}r_{\RRR}(p'')\;<\;r_{\RRR}(p'_j)\;\le\;r_{\RRR}(p'').
\ee
Then
\be
\begin{aligned}
\frac{\Vol{B}(p'',\,\frac{3l}{k}r_{\RRR}(p''))}{\Vol{B}(p'_j,\,\frac{1}{2k}r_{\RRR}(p'_i))} 
&\;\le\;
\frac{\Vol{B}\left(p'',\,\left(\frac{3l}{k}\frac{k+l}{k-l}\right)r_{\RRR}(p'_j)\right)}{\Vol{B}_{p'_j}(\frac{1}{2k}r_{\RRR}(p'_j))} \\
&\;\le\; \frac{\Vol{B}\left(p'_j,\,\left(\frac{3l}{k}\frac{k+l}{k-l}\right)r_{\RRR}(p'_j)+\dist(p'',p'_j)\right)}{\Vol{B}_{p'_j}(\frac{1}{2k}r_{\RRR}(p'_j))} \\
&\;=\;\frac{\Vol{B}\left(p'_j,\,\frac{l}{k}\left(1\,+\,4\frac{k+l}{k-l}\right)r_{\RRR}(p'_j)\right)}{\Vol{B}(p'_j,\,\frac{1}{2k}r_{\RRR}(p'_i))}.
\end{aligned}
\ee
Then as long as 
\be
{B}\left(p'_j,\,\frac{l}{k}\left(1\,+4\frac{k+l}{k-l}\right)r_{\RRR}(p'_j)\right) \subset {B}(p'_j,\,r_{\RRR}(p'_j)),
\ee
(which holds if $l\le\frac17k$, or more precisely if $l\le\frac{1}{3+\sqrt{12}}k$), we can use volume comparison to obtain
\be
\begin{aligned}
\frac{\Vol{B}(p'',\,\frac{3l}{k}r_{\RRR}(p''))}{\Vol{B}(p'_j,\,\frac{1}{2k}r_{\RRR}(p'_i))}
&\;\le\; C(n)\frac{\left(\frac{l}{k}\left(1\,+\,4\frac{k+l}{k-l}\right)r_{\RRR}(p'_j)\right)^n}{(\frac{1}{2k}r_{\RRR}(p'_j))^n} \\
&\;=\; C(n)\,2^n\,l^n\left(1\,+\,4\frac{k+l}{k-l}\right)^n \\
&\;\le\; C'(n)\,k^n
\end{aligned}
\ee
which, combined with (\ref{IneqMultBound}), gives the conclusion.
\qed

\begin{corollary}[Gromov-style packing lemma with a cutoff] \label{CorGromovCoverII}
Let $\Omega\subseteq{N}^n$ be a domain in a Riemannian $n$-manifold, and let $\{p_i\}$ be a maximally $\frac1k{r}_{\RRR}$-separated subset of $\Omega$ where $k>1$.
Assume $\frac{k}{k-1}<l\le\frac17k$.
Then
\be
\Omega \;\subseteq\; \bigcup_{i=1}^N\;B\left(p_i,\,\frac{l}{k}r_{\RRR}^s(p_i)\right) \;\subset\;\Omega^{\mathcal{R},\frac{l}{k}s}
\ee
and the multiplicity of this cover is uniformly bounded by $C(n)\cdot{k}^n$.
\end{corollary}
{\it Proof}.
From the previous lemma we have a maximally $r_{\mathcal{R}}$-separated collection of points $P=\{p_i\}_{i=1}^N$ (possibly $N=\infty$) so that the union of $B(p_i,\frac{l}{k}r_{\mathcal{R}}(p_i))$ covers $\Omega$.
Let $P^{\mathcal{R}}\subseteq{P}$ be the set
\be
P^{\mathcal{R}}\;=\;\left\{\;p\in{P}\;\big|\;r_{\mathcal{R}}^s(p)\,<\,s\;\right\}\;=\;\left\{p_i^{\mathcal{R}}\right\}_{i=1}^{N^{\mathcal{R}}}.
\ee
Let
\be
\Omega'\;=\bigcup_{i=1}^{\;\;\;N^{\mathcal{R}}}B\left(p_i,\,\frac{l}{k}r^s_{\mathcal{R}}(p_i)\right).
\ee
Now if $q$ is any point in $\Omega'\setminus\Omega$, then $p$ is a distance of at most $\frac{l}{k}r_{\mathcal{R}}(p_j)<\frac17r_{\mathcal{R}}(p_j)$ from some point $p_j$ where $r_{\mathcal{R}}(p_j)\ge{s}$; therefore $r_{\mathcal{R}}(q)\ge\frac67s$.
Therefore the enlarged set $\left(\Omega\setminus\Omega'\right)^{\left(\frac67s\right)}$ has $|\Riem|\le\frac{49}{36}s^{-2}$.

Scaling the metric by $(\frac67s)^{-1}$, we have $|\Riem|\le1$ on $\left(\Omega\setminus\Omega'\right)^{(1)}$.
Now the usual Gromov packing lemma provides a covering of $\Omega\setminus\Omega'$ by a collection of balls $\{B(q'_i,\frac{l}{k})\}_{i=1}^{N^s}$ where $q_i\in\Omega\setminus\Omega'$, and there is a multiplicity bound for this cover.
Returning to the original scale, we have the cover $\{B(q_i,\frac67\frac{l}{k}s)\}_{i=1}^{N^s}$, and since $\frac67\frac{l}{k}s\le\frac{l}{k}r^{s}_{\mathcal{R}}(q_i)$, we have that
\be
P^s\;\triangleq\;\left\{B(q_i,\frac{l}{k}r^s_{\mathcal{R}}(p))\,\right\}
\ee
covers $\Omega\setminus\Omega'$.
Finally we note that the multiplicity of the overall covering $P^{\mathcal{R}}\cup{P}^s$ is bounded by the sum of the multiplicities of the individual covers, we have the multiplicity controlled by $C(n)k^n+C(n)$ which is controlled by $C(n)k^n$, as $k\ge8$.
\qed

\subsection{The integration lemma} \label{SubSecIntegrationLEmma}

Given $\mu,s>0$, we modify the definition of $\Omega^{\mathcal{R},s}$ to include another scale factor:
\be
\Omega^{\mu\mathcal{R},s}
\;=\;
\left\{\;p\in{N}\;\big|\;p\,\in\,B\left(p,\,\mu\,r_{\mathcal{R}}^s(p)\right)\;\right\}.
\ee
Note $\Omega^{\mu\mathcal{R},s}\subseteq\Omega^{\mathcal{R},\mu{s}}$.
We shall prove the following improved proposition:
\begin{proposition} \label{PropIntegrationLemmaII}
Given $K<\infty$, $\mu\in(0,1]$, and $k>0$, there exist numbers $\epsilon_0=\epsilon(K)>0$, $C=C(\mu,K,k)<\infty$, and a universal constant $m<\infty$, so that $s<{K}\Lambda^{-1}$ and $\int_{\Omega^{(2s)}}|\Riem|^2\;\le\;\epsilon_0$ together imply
\be
\int_{\Omega}\left(r^s_{\mathcal{R}}\right)^{-k}
\;\le\;m\,{s}^{-k}\,\Vol\,\Omega^{(\mu{s})}
\;+\;{C}\;\int_{\Omega^{\mu\mathcal{R},s}}|\Riem|^{\frac{k}{2}}. \label{IneqIntLemmaPrecise}
\ee
Further, if $\int_{\Omega^{2\mathcal{R}}}|\Riem|^2<\epsilon_0$ and either $\Ric\ge0$ or else $r_{\mathcal{R}}(p)<K\Lambda^{-1}$ for all $p\in\Omega$, we have
\be
\int_{\Omega}\left(r_{\mathcal{R}}\right)^{-k}
\;\le\;{C}\;\int_{\Omega^{\mu\mathcal{R}}}|\Riem|^{\frac{k}{2}}
\ee
(where $C$ does not depend on $K$ in the case $\Ric\ge0$).
\end{proposition}
{\it Proof}. 
Consider a domain $\Omega$ and a number $s>0$, so that $\int_{\Omega^{(\frac\mu4{s})}}|\Riem|^2<\epsilon_0$.
Using the covering lemma, specifically Corollary \ref{CorGromovCoverII}, with $k=8\mu^{-1}$ and $l=8/7$, we have a set of points $P=\{p_i\}_{i=1}^N\subseteq\Omega$ so that
\be
\Omega\;\subset\;\bigcup_{i=1}^N{B}(p_i,\,\mu\frac17r_{\mathcal{R}}^s) \;\subset\Omega^{\mu\frac17\mathcal{R},\mu\frac17s}
\ee
and so that the multiplicity is less than some dimensional constant $m$.
Now we divide the points $p_i$ into two sets, 
\be
\begin{aligned}
&P^{\mathcal{R}}\;=\;\left\{\;p\,\in\,{P}\;\big|\;r_{\mathcal{R}}^s\,<\,s\;\right\}\;=\;\{p_i^R\}_{i=1}^{N^{\mathcal{R}}} \\
&P^s\;=\;\left\{\;p\,\in\,{P}\;\big|\;r_{\mathcal{R}}^s\,=\,s\;\right\}\;=\;\{p_i^s\}_{i=1}^{N^s}
\end{aligned}
\ee
Due to Lipschitz control on $r_{\mathcal{R}}^s$ we have that $r_{\mathcal{R}}^s\ge(1-\mu\frac17)r_{\mathcal{R}}^s(p_i)$ on $B(p_i,\mu\frac17r_{\mathcal{R}}^s(p_i))$, so
\be
\Xint{\;-}_{B(p_i,\,\mu\frac17{r}_{\mathcal{R}}^s(p_i))}\left(r_{\mathcal{R}}^s\right)^{-k}\;\le\;\frac67\left(r_{\mathcal{R}}^s(p_i)\right)^{-k}
\ee
If $p_i\in{P}$ is in fact in $P^{\mathcal{R}}$ then the second Harnack inequality, Theorem \ref{ThmHarnackII}, now gives
\be
\begin{aligned}
&\int_{B\left(p^{\mathcal{R}}_i,\,\frac{\mu}{7}r_{\mathcal{R}}^s(p^{\mathcal{R}}_i)\right)}\left(r_{\mathcal{R}}^s\right)^{-k}\;\le\;C\int_{B\left(p_i^{\mathcal{R}},\,\frac{\mu}{14}r_{\mathcal{R}}(p^{\mathcal{R}}_i)\right)}|\Riem|^{\frac{k}{2}}
\end{aligned}
\ee
where $C=C(\mu,K,k)$, provided $\epsilon=\epsilon(K)$ is sufficiently small.
Recalling that the half-radius balls $B\left(p_i^{\mathcal{R}},\,\frac{\mu}{14}r_{\mathcal{R}}(p^{\mathcal{R}}_i)\right)$ are disjoint, we have
\be
\begin{aligned}
\int_{\Omega}\left(r_{\mathcal{R}}^s\right)^{-k}
&\;\le\;\sum_{i=1}^N\int_{B\left(p_i,\,\frac{\mu}{7}r_{\mathcal{R}}^s(p_i)\right)}\left(r_{\mathcal{R}}^s\right)^{-k} \\
&\;=\;
\sum_{i=1}^{N^s}\int_{B\left(p^s_i,\,\frac\mu8r_{\mathcal{R}}^s(p_i^s)\right)}\left(r_{\mathcal{R}}^s\right)^{-k} \;+\;
\sum_{i=1}^{N^{\mathcal{R}}}\int_{B\left(p^{\mathcal{R}}_i,\,\frac\mu8r_{\mathcal{R}}^s(p_i^{\mathcal{R}})\right)}\left(r_{\mathcal{R}}^s\right)^{-k} \\
&\;\le\;\sum_{i=1}^{N^s}s^{-k}\,\Vol\,B\left(p_i^s,\,\frac{\mu}{7}s\right) \;+\;
C\sum_{i=1}^{N^{\mathcal{R}}}\int_{B\left(p_i^R,\,\frac{\mu}{14}r_{\mathcal{R}}^s(p_i^{\mathcal{R}})\right)}|\Riem|^{\frac{k}{2}} \\
&\;\le\;{m}\,s^{-k}\,\Vol\,\Omega^{(\frac{\mu}{7}s)} \;+\;
C\,m\,\int_{\Omega^{\frac{\mu}{14}\mathcal{R},s}}|\Riem|^{k/2}.
\end{aligned}
\ee
where $C=C(\mu,K,k)$ and $m$ is the universal multiplicity bound given by Corollary \ref{CorGromovCoverII}.
Now if $s=\infty$, we have that $P^s$ is empty, and we have just
\be
\begin{aligned}
\int_{\Omega}\left(r_{\mathcal{R}}\right)^{-k}
&\;\le\;C\,m\,\int_{\Omega^{\frac{\mu}{14}\mathcal{R}}}|\Riem|^{k/2}.
\end{aligned}
\ee

\section{Proof of the epsilon regularity theorem} \label{PfEpsReg}

\subsection{Construction of the Transgression} \label{SubSecTransgression}

If $\int_{B(p,r)}|\Riem|^2<\epsilon_0$, for some sufficiently small $\epsilon_0$, and if $r<K\Lambda^{-1}$, we shall construct a transgression (on a slightly smaller sub-ball) for any characteristic 4-form $\mathcal{P}$.
By Lemma \ref{LemmaChTiForcedCollpasing}, $B(p,\frac12r)$ is $\tau$-collapsed on the scale $r_{\mathcal{R}}(p)$, and so by Theorem \ref{ThmCollapseBigThm} has an N-structure $\mathfrak{N}=(\Omega,\mathcal{N},\iota)$, where $\Omega$ is some saturation of $B(p,\frac12r)$.

{\it \underline{Step I}: Modification of the Levi-Civita connection of an invariant metric}

Let $g'$ be the invariant metric on $\Omega$, and $\nabla'$ its Levi-Civita connection; it is $C^{k,\alpha}$-close to $g$ in the scale $r_{\mathcal{R}}^s(p)$.
Let $(U_l,\mathfrak{F}_l,\iota_l)$ be an atlas, which we can assume is $C$-regular.
For the moment, let us pass to the universal covers $\widetilde{U_l}$, where the sheafs $\mathcal{F}_l$ are represented by Killing fields $v_1,\dots,v_k$ defined on $\widetilde{U}_l$.
We may assume the fields $v_1,\dots,v_k$ are mutually orthogonal at a point $p$, and since they commute, we may assume they are globally orthogonal.
Because the F-structure is polarizeable, we cam assume that the $\{v_i\}$ are all nowhere-zero on their domains of definition.
Let $K^{(l)}$ be the transformation-valued $1$-form
\begin{eqnarray}
K^{(l)} &=& \sum_i|v_i|^{-2}\left(v_i\right)_\flat\,\nabla{v}_i. \label{DefKl}
\end{eqnarray}
Because $\nabla{v}_i$ is anti-symmetric, we may consider $K^{(l)}$ to be an $\mathfrak{so}(4)$-valued 1-form on $\widetilde{U_l}$.
Because $K^{(l)}$ is invariant under any (constant-coefficient) orthogonal transformation of the $v_i$, it is invariant under the deck transformations.
Thus it passes back to an $\mathfrak{so}(4)$-valued 1-form on $U_l$.

Let $\{\varphi_l\}$ be the partition of unity subordinate to $\{U_l\}$ we can assume that $|\nabla\varphi_l|<Cr_{\mathcal{R}}$.
Define the $\mathfrak{so}(4)$-valued 1-form $K$ by
\begin{eqnarray}
K &=& \sum_l\;\varphi_l\,{K}^{(l)}.
\end{eqnarray}
Define a connection $\widetilde{\nabla}=\nabla'-K$.
Its connection 1-form is $\widetilde{A}=A'-K$, and its curvature 2-form is $\widetilde{F}=d\widetilde{A}+\frac12[\widetilde{A},\widetilde{A}]$.

{\it \underline{Step II}: Transgression of characteristic 4-forms of the invariant metric}

Any is any $SO(4)$-invariant polynomial $\mathcal{P}$ of degree 2 on $\mathfrak{so}(4)$, leads to a characteristic 4-form on a 4-manifold by plugging in the curvature 2-form associated to any connection.
First we consider the 4-forms $\mathcal{P}(\widetilde{F},\widetilde{F})$ associated to the modified connections $\widetilde{A}$.

If $\{v_j\}$ are the orthogonal set of Killing fields on $\Omega^{(l)}$, then $i_{v_j}K=\nabla'{v}_j$.
This means that if a locally defined Killing field $v$ near a point $p$ is common to all the structures $\mathfrak{F}_l$ on all $\Omega^{(l)}$ that contain $p$, then $\widetilde{\nabla}v=\nabla'v-\sum_l{i}_{v}\varphi_l{K}^{(l)}=0$, so $i_v\widetilde{A}=0$.
We have then
\be
\begin{aligned}
i_v\widetilde{F}
&\;=\;i_vd\widetilde{A}\,+\,\frac12i_v[\widetilde{A},\,\widetilde{A}] \\
&\;=\;\mathcal{L}_v\widetilde{A}\,-\,di_v\widetilde{A}\,+\,[i_v\widetilde{A},\,\widetilde{A}]
\;=\;0
\end{aligned}
\ee
where $\mathcal{L}$ is the Lie derivative.
Now from property ({\it{c}}) of Theorem \ref{ThmCollapseBigThm} we have that the stalks on any overlaps are ordered by containment, so there is indeed a common locally defined Killing field around every point $p\in\Omega$.
Therefore $\widetilde{F}$ has at least one null vector at every point.
Therefore $\mathcal{P}(\widetilde{F},\,\widetilde{F})$ has a null vector, and therefore $\mathcal{P}(\widetilde{F},\,\widetilde{F})=0$.

Going through the usual interpolation procedure, set $A'_t=A'(1-t)+\widetilde{A}t=A'-Kt$, so $F_t=dA_t+\frac12[A_t,A_t] =F'\,-\,tD'K\,+\,\frac12t^2[K,K]$.
Then
\be
0\;=\;\mathcal{P}(\widetilde{F},\widetilde{F})
\;=\;\mathcal{P}(F',F')\;+\;\frac12d\int_0^1\mathcal{P}(K,\,F_t')\,dt.
\ee
To come to an explicit estimate, set $\mathcal{TP}_1=\frac12\int_0^1\mathcal{P}(K,\,F_t')\,dt$.
Because
\be
D'K^{(l)}\;=\;\sum_id(v_i)_\flat\otimes\left(|v_i|^{-2}\nabla{v}_i\right)\,+\,(v_i)_\flat\wedge{D}'\left(|v_i|^{-2}\nabla{v}_i\right)
\ee
we have that
\be
\begin{aligned}
\mathcal{P}(K^{(l)},D'K^{(l)})
\;=\;\sum_{i,j}|v_i|^{-2}|v_j|^{-2}(v_i)_\flat\wedge{d}(v_j)_\flat\;\cdot\;\mathcal{P}(\nabla{v}_i,\,\nabla{v}_j)
\end{aligned}
\ee
so we can estimate, at a given point $p$, that
\be
\begin{aligned}
\left|\mathcal{TP}_1\right|
&\;\le\;
C\left(|\Riem|\;\cdot\;{max}_i\frac{|\nabla{v_i}|}{|v_i|}\;+\;\left({max}_i\frac{|\nabla{v_i}|}{|v_i|} \;+\;\,{max}_j|d\varphi_j|\right)\;{max}_i\frac{|\nabla{v_i}|^2}{|v_i|^2}\right) \\
&\;\le\;C\,r_{\mathcal{R}}^r(p)^{-3}
\end{aligned}
\ee

{\it \underline{Step III}: Construction of the transgression for the original metric}

By ({\it{e}}) of Theorem \ref{ThmCollapseBigThm}, if $A$ is the connection 1-form in the original metric, then setting $L=A-A'$, we have $|L|$ is $C^{k,\alpha}$-close to zero at $p$, after the metric is scaled so that $r_{\mathcal{R}}^r(p)=1$.
Therefore $|L|<Cr_{\mathcal{R}}^r(p)^{-1}$.
Setting $F_t=dA_t+\frac12[A_t,A_t]=F-tDL+\frac12t^2[L,L]$ so
\be
\mathcal{TP}_2\;=\;\frac12\int_0^1\mathcal{P}(K,\,F_t)\,dt
\ee
is also $C^{k,\alpha}$-close to zero if the metric is scaled so $r_{\mathcal{R}}(p)=1$, so that $|\mathcal{TP}_2|\le{C}r_{\mathcal{R}}^r(p)^{-3}$.
Setting $\mathcal{TP}=\mathcal{TP}_1+\mathcal{TP}_2$, we have
\be
\mathcal{P}(F,\,F)\;+\;d\mathcal{TP}\;=\;0 \label{EqnTransgressionEqn}
\ee
and $|\mathcal{TP}|\le{C}\left(r_{\mathcal{R}}^s(p)\right)^{-3}$ for any $s<r<K\Lambda^{-1}$.

\subsection{Use of the Transgressions}

We have two canonical $SO(4)$ characteristic forms, $\mathcal{P}_\chi$ and $\mathcal{P}_\tau$ corresponding to the Euler form and the Pontrjagin form (three times the signature form $\mathcal{P}_\tau$), where
\be
\begin{aligned}
\mathcal{P}_\chi
&\;=\;\frac{1}{8\pi^2}\left(\frac{1}{24}R^2\,-\,\frac12\big|\cRic\big|^2\,+\,\frac14\left|W^{+}\right|^2\,+\,\frac14\left|W^{-}\right|^2\right) \\
\mathcal{P}_\tau
&\;=\;\frac{1}{12\pi^2}\left(\frac14\left|W^{+}\right|^2\,-\,\frac14\left|W^{-}\right|^2\right).
\end{aligned}
\ee
(the unusual factor of $\frac14$ comes from our use of the tensor norm, instead of the more common $\bigwedge^2$-norm induced by the Hodge star).
With $|\Riem|^2=\frac16R^2+2|\cRic|^2+|W|^2$ we have
\be
|\Riem|^2\;=\;\frac13R^2\,+\,4\left|W^+\right|^2\,-\,32\pi^2\left(\mathcal{P}_\chi\,+\,3\mathcal{P}_\tau\right)
\ee
This is most usable if $|W^+|$ has apriori control, as is the case when the metric is half conformally flat whence $W^+=0$, or K\"ahler whence $|W^+|^2=\frac16R^2$.
Writing $\mathcal{P}=32\pi^2(\mathcal{P}_\chi+\mathcal{P}_\tau)$ the previous section gives a 3-form $\mathcal{TP}$ with $|\mathcal{TP}|<C\left(r_{\mathcal{R}}^s\right)^{-3}$ when $s\le{K}\Lambda^{-1}$, as long as we remain within $\Omega$.
We have
\be
|\Riem|^2\;=\;\frac13{R}^2\,+\,4\left|W^+\right|^2\,-\,d\mathcal{TP}
\ee
where $\alpha=1$ in the K\"ahler case and $\alpha=\frac13$ in the half-conformally flat case.
If $\varphi$ is a $C_c^\infty$ function with support in $\Omega$ then
\be
\begin{aligned}
\int\varphi|\Riem|^2
&\;=\;\int\varphi\left(\frac13{R}^2+4\left|W^+\right|^2\right)\,+\,\int{d}\varphi\wedge\mathcal{TP} \\
&\;\le\;\int\varphi\left(\frac13{R}^2+4\left|W^+\right|^2\right)\,+\,C\int|\nabla\varphi|\left(r_{\mathcal{R}}^s\right)^{-3} \label{IneqChTiStartI}
\end{aligned}
\ee
where $C=C(K)$ and $s\le{K}\Lambda^{-1}$.

After picking a basepoint $q$, we immediately suppress it and write $B(\rho)$ for $B(q,\rho)$ and denote annuli by $A(\rho',\rho)=B(q,\rho)\setminus{B}(q,\rho')$.
To avoid using ``$Vol$'' frequently, we use the symbol $|\cdot|$ to indicate Hausdorff $4$-volume.
From volume comparison we have
\be
\frac{\left|A(\rho-\gamma,\rho+\beta)\right|}{\left|B(\rho)\right|}
\;\le\;\frac{\left|B(\rho+\beta)\right|}{\left|B(\rho)\right|}
\;\le\;C\,\left(1+\frac\beta\rho\right)^4\,\cosh^3(\Lambda\beta)  \label{IneqChTiStartII}
\ee
where we can take for instance $C=\frac98$.
The integration lemma (Proposition \ref{PropIntegrationLemmaII}) and the inequalities in (\ref{IneqChTiStartI}) and (\ref{IneqChTiStartII}) allow us to carry out the Cheeger-Tian iteration process.

\subsection{Proof of ({\it{i}})}

Here we assume $r$ is too big compared to Ricci curvature: $r>K\Lambda^{-1}$.
For technical convenience we will assume $K>100$.
Pick the basepoint $q$ to lie within $B(p,\frac12r)$, which contains $B(p,50\Lambda^{-1})$.

Select $\rho=\Lambda^{-1}$ and choose any $\mu\le1$.
Let $\varphi$ be a cutoff function with support in $B(\rho+\frac12\mu\Lambda^{-1})$, so that $|\nabla\varphi|<4\mu\Lambda^{-1}$ and $\varphi\equiv1$ in $B(\rho)$.
By (\ref{IneqChTiStartI})
\be
\begin{aligned}
\int_{B(\rho)}|\Riem|^2
&\;\le\;\alpha\int_{B(\rho+\frac12\mu\Lambda^{-1})}\left(\frac13{R}^2+4\left|W^+\right|^2\right) \\
& \quad\quad \,+\,C\left(\mu\Lambda^{-1}\right)^{-1}\int_{A(\rho,\rho+\frac12\mu\Lambda^{-1})}\left(r_{\mathcal{R}}^{\frac12\mu\Lambda^{-1}}\right)^{-3}
\end{aligned}
\ee
The integration lemma Proposition \ref{PropIntegrationLemmaII} gives us
\be
\begin{aligned}
\int_{A(\rho,\rho+\frac12\mu\Lambda^{-1})}\left(r_{\mathcal{R}}^{\frac12\mu\Lambda^{-1}}\right)^{-3}
&\;\le\;C\left(\mu\Lambda^{-1}\right)^{-3}\left|A(\rho-\frac12\mu\Lambda^{-1},\rho+\mu\Lambda^{-1})\right| \\
&\quad\quad\,+\,C\int_{A(\rho-\frac12\mu\Lambda^{-1},\rho+\mu\Lambda^{-1})}|\Riem|^{\frac32} \\
\end{aligned}
\ee
where $C$ is universal (does not depend on $K$).
Volume comparison and H\"older's inequality allows the average curvature estimate
\be
\begin{aligned}
\Xint{\;-}_{B(\rho)}|\Riem|^2
&\;\le\;\Xint{\;-}_{B(\rho+\frac12\mu\Lambda^{-1})}\left(\frac13{R}^2+4\left|W^+\right|^2\right)\,+\,
C'\left(\mu\Lambda^{-1}\right)^{-4}\cosh^3(\mu) \\
&\quad\quad\,+\,C\left(\mu\Lambda^{-1}\right)^{-1}\cosh^3(\mu)\Xint{\;-}_{B(\rho+\mu\Lambda^{-1})}|\Riem|^{\frac32} \\
&\;\le\;\Xint{\;-}_{B(\rho)}\left(\frac13{R}^2+4\left|W^+\right|^2\right)\,+\,C\left(\mu\Lambda^{-1}\right)^{-4}\left(\cosh^3(\mu)+\cosh^{12}(\mu)\right) \\
&\quad\quad\quad\quad\quad\quad \,+\,\frac34\Xint{\;-}_{B(\rho+\mu\Lambda^{-1})}|\Riem|^2.
\end{aligned}
\ee
Now we iterate.
Choose any $\rho_0=\rho=\Lambda^{-1}$ and $\rho_{i+1}=\rho_i+\mu_i\Lambda^{-1}$ for some sequence $\{\mu_i\}$ with, say, $\sum_i\mu_i<25$.
Then we retain $\int_{B(2\rho_i)}|\Riem|^2<\epsilon_0$, so for any $T\in\mathbb{N}$ we have
\be
\begin{aligned}
\Xint{\;-}_{B(\rho_0)}|\Riem|^2
&\;\le\;C\Xint{\;-}_{B(\rho_{T})}\left(\frac13{R}^2+4\left|W^+\right|^2\right)\,+\,
C\Lambda^4\sum_{i=0}^{T-1}\left(\frac34\right)^i\mu_i^{-4} \\
&\quad\quad\quad\quad\quad\quad
\,+\,\left(\frac34\right)^T\Xint{\;-}_{B(\rho_{T})}|\Riem|^2.
\end{aligned}
\ee
Setting $\mu_i=\left(\frac{33}{40}\right)^{\frac{i}{4}}$ so $\frac{3^i}{4^i}\mu_i^{-4}=\left(\frac{10}{11}\right)^i$, we have $\sum\frac{3^i}{4^i}\mu_i^{-4}=11$ and $\lim_i\rho_i\approx20.3\Lambda^{-1}$.
Taking $T\rightarrow\infty$ we obtain
\be
\begin{aligned}
\Xint{\;-}_{B(\rho_0)}|\Riem|^2 &\;\le\;C\Xint{\;-}_{B(21\Lambda^{-1})}\left(\frac13{R}^2+4\left|W^+\right|^2\right)\,+\,C\Lambda^4.
\end{aligned}
\ee
Now if $r_{\mathcal{R}}\ge(2C)^{-\frac14}\Lambda^{-1}$ for any point in $B(\Lambda^{-1})$, then the first Harnack inequality Theorem \ref{ThmHarnackI} implies $r_{\mathcal{R}}(q)\ge{C}'\Lambda^{-1}$ on $B(\Lambda^{-1})$ for some $C'$.
On the other hand if $r_{\mathcal{R}}<(2C)^{-\frac14}\Lambda^{-1}$ for all point in $B(\Lambda^{-1})$, then the second inequality of the integration lemma along with volume comparison implies
\be
\begin{aligned}
\Xint{\;-}_{B(\Lambda^{-1})}\left(r_{\mathcal{R}}\right)^{-4}
&\;\le\;C\Xint{\;-}_{B(\Lambda^{-1})}|\Riem|^2 \\
&\;\le\;C\Xint{\;-}_{B(21\Lambda^{-1})}\left(\frac13{R}^2+4\left|W^+\right|^2\right)\,+\,C\Lambda^4
\end{aligned}
\ee
which implies that $\left(r_{\mathcal{R}}(q')\right)^{-4}$ is bounded for some $q'\in{B}(\Lambda^{-1})$.
Therefore
\be
2C\Lambda^4\;<\;
r_{\mathcal{R}}(q)^{-4}\;<\;C\Xint{\;-}_{B(21\Lambda^{-1})}\left(\frac13{R}^2+4\left|W^+\right|^2\right)\,+\,C\Lambda^4
\ee
so we directly obtain $\Xint{\;-}_{B(21\Lambda^{-1})}\left(\frac13{R}^2+4\left|W^+\right|^2\right)>\Lambda^4$, and the conclusion follows.

\subsection{Proof of ({\it{ii}})}

Now assume $r\le{K}\Lambda^{-1}$.
Choose $\mu$, $\rho$ so $\mu<\frac{1}{100}r$ and choose $\rho<\frac{1}{100}r$.
Repeating the procedure above, we obtain the estimate
\be
\begin{aligned}
\Xint{\;-}_{B(\rho)}|\Riem|^2
&\;\le\;C\mu^{-4}
\,+\,\Xint{\;-}_{B(\rho+\mu)}\left(\frac13{R}^2+4\left|W^+\right|^2\right)
\,+\,\frac34\Xint{\;-}_{B(\rho+\mu)}|\Riem|^2.
\end{aligned}
\ee
where now $C$ depends on $K$, via $\cosh(\mu\Lambda)\le\cosh(K/8)$.
Iterating again gives
\be
\begin{aligned}
\Xint{\;-}_{B(\rho_0)}|\Riem|^2
&\;\le\;C\Xint{\;-}_{B(\rho_{T})}\left(\frac13{R}^2+4\left|W^+\right|^2\right) \\
&\quad\quad\,+\,C\sum_{i=0}^{T-1}\left(\frac34\right)^i\mu_i^{-4}
\,+\,\left(\frac34\right)^T\Xint{\;-}_{B(\rho_{T})}|\Riem|^2 \\
\Xint{\;-}_{B(r/100)}|\Riem|^2
&\;\le\;C\Xint{\;-}_{B(r)}\left(\frac13{R}^2+4\left|W^+\right|^2\right)\,+\,
C\,r^{-4}
\end{aligned}
\ee
where we have chosen $\rho_{i+1}=\rho_i+\mu_i$, $\rho_0=r/100$, and $\mu_i=r\left(\frac{33}{40}\right)^{\frac{i}{4}}/25$.

Reasoning as above, if some $q\in{B}(p,\frac12r)$ exists with $r_{\mathcal{R}}(q)>(2C)^{-1/4}r$, then the Harnack inequality Theorem \ref{ThmHarnackI} implies $r_{\mathcal{R}}>C'r$ on all of $B(p,\frac12r)$ for some $C'=C'(K)$ which implies $|\Riem|<C'r^{-2}$ on $B(p,\frac12r)$.
If not, meaning $r_{\mathcal{R}}(q)<<r$ for all $q\in{B}(p,\frac12r)$, then the integration lemma gives
\be
2Cr^{-4}\;<\;r_{\mathcal{R}}(q)^{-4}\;<\;C\Xint{\;-}_{B(r)}\left(\frac13{R}^2+4\left|W^+\right|^2\right)\,+\,Cr^{-4}
\ee
so that $\Xint{\;-}_{B(r)}R^2>Cr^{-4}$ and $sup\,|\Riem|^2<sup\,(r_{\mathcal{R}})^{-4}<C\Xint{\;-}_{B(r)}\left(\frac13{R}^2+4\left|W^+\right|^2\right)$.

\subsection{Proof of ({\it{iii}})} \label{SubSectionPfiii}

This only has additional meaning if
\be
\Xint{\;-}_{B(p,r)}|\Riem|^2\;<\;\delta\,r^{-4}.
\ee
for some $\delta<<1$; indeed choose $\delta=C^{-1}$.
First, we rule out that every point $q\in{B}(p,\frac12r)$ has $r_{\mathcal{R}}(q)<\delta_0r$, where $\delta_0$ is the constant provided in the Harnack ineequality, Theorem \ref{ThmHarnackI}.
The integration lemma provides a constant $C$ so that
\be
\Xint{\;-}_{B(q,\frac14r)}\left(r_{\mathcal{R}}\right)^{-4}
\;<\;C\Xint{\;-}_{B(q,\frac12r)}|\Riem|^2
\;<\;C\cdot\delta\,r^{-4}\;=\;r^{-4}.
\ee
Thus $r_{\mathcal{R}}>(C\delta)^{-\frac14}r=r$ somewhere on $B(q,\frac12r)$, and so by the Harnack inequality $r_{\mathcal{R}}>\delta_0r$ on $B(p,\frac12r)$, an absudity.

Now suppose $|\Riem(q)|^2>>\int_{B(p,r)}|\Riem|^2$ at some $q\in{B}(p,\frac12r)$.
Because $r_{\mathcal{R}}(p)>\delta_0r$, the exponential map has no conjugate points on a ball of radius, say, $\delta_0(1+\pi/10)r$.
Lifting to the (contractible) ball of the same radius in the tangent space at $q$, we can apply the standard elliptic thoery to obtain $|\Riem|^2<C'\Xint{\;-}_{B(p,\delta_0r)}|\Riem|^2$ for some universal $C'$.
Passing back down to the manifold we retain this inequality for a slightly different $C'$.
Then volume comparison gives
\be
|\Riem(q)|^2\;\le\;C'\Xint{\;-}_{B(q,\delta_0r)}|\Riem|^2\;\le\;C\Xint{\;-}_{B(p,r)}|\Riem|^2
\ee
where $C$ now depends also on $K$.

\subsection{Proof of Corollary \ref{CorHarnack}}

First, if $r_{\mathcal{R}}(p)$ is large compared to $r$ at even one point, then the Harnack inequality (Theorem \ref{ThmHarnackI}) says it is comparable to $r$ at all points in the half-radius ball, so the conclusion holds in this case.

On the other hand if $r_{\mathcal{R}}$ is small compared to $r$ everywhere, are in the situation of the proof of ({\it iii}) above, and we obtain
\be
|\Riem(q)|^2\;\le\;C\Xint{\;-}_{B(p,r)}|\Riem|^2
\ee
for all $q\in{B}(p,\frac12r)$, so that $\sup_{B(p,\frac14r)}\left(r_{\mathcal{R}}\right)^{-4}$ satisfies the same inequality.
We have trivially that
\be
\int_{B(p,r)}|\Riem|^2\;\le\;\int_{B(p,r)}\left(r_{\mathcal{R}}\right)^{-4},
\ee
which completes the assertion.

\end{document}